\documentclass[12pt,reqno]{article}

\usepackage[usenames]{color}
\usepackage{amssymb}
\usepackage{amsmath}
\usepackage{amsthm}
\usepackage{amsfonts}
\usepackage{amscd}
\usepackage{graphicx}

\usepackage[colorlinks=true,
linkcolor=webgreen,
filecolor=webbrown,
citecolor=webgreen]{hyperref}

\definecolor{webgreen}{rgb}{0,.5,0}
\definecolor{webbrown}{rgb}{.6,0,0}

\usepackage{color}
\usepackage{fullpage}
\usepackage{float}

\usepackage{breakurl}

\newcommand{\Mod}[1]{\ \mathrm{mod}\ #1}
\begin{document}

\theoremstyle{plain}
\newtheorem{theorem}{Theorem}
\newtheorem{corollary}[theorem]{Corollary}
\newtheorem{lemma}[theorem]{Lemma}
\newtheorem{proposition}[theorem]{Proposition}
\theoremstyle{definition}
\newtheorem{definition}[theorem]{Definition}
\newtheorem{example}[theorem]{Example}
\newtheorem{conjecture}[theorem]{Conjecture}

\theoremstyle{remark}
\newtheorem{remark}[theorem]{Remark}

\begin{center}
\vskip 1cm{\LARGE\bf A Twist of a Ramanujan Identity
\vskip 1cm}
\large

H\`ung Vi\d{\^e}t Chu\\
Department of Mathematics\\
University of Illinois at Urbana-Champaign \\
Champaign, IL 61820\\
USA \\
\href{mailto:hungchu2@illinois.edu}{\tt hungchu2@illinois.edu} \\
\ \\
L\^an Kh\'anh Chu\\
Research and Consultancy Department\\
Banking Academy of Vietnam\\
Hanoi 10000\\
Vietnam \\
\href{mailto:lanck@hvnh.edu.vn}{\tt lanck@hvnh.edu.vn} \\
\end{center}
\vskip .2 in

\begin{abstract}
Ramanujan wrote the following identity
\begin{align*}
\sqrt{2 \left(1 - \frac{1}{3^2}\right)
\left(1 - \frac{1}{7^2}\right)
\left(1 - \frac{1}{11^2}\right)
\left(1 - \frac{1}{19^2}\right)}
\  =
\
 \left(1 + \frac{1}{7}\right) \left(1 + \frac{1}{11}\right) \left(1 + \frac{1}{19}\right), 
\end{align*} on which Berndt asked "Is this an isolated result, or are there other identities of this type?" Reb\'ak provided formulas that generate many similar identities and believed that the curious identity is related to the reciprocal of the Landau-Ramanujan constant. In a previous work, Chu examined necessary and sufficient conditions for the integers in the identity and proved that there are only finitely many such identities. In this note, we twist the identity to have infinitely many Ramanujan-type identities.
\end{abstract}

\section{Introduction}
In ranking mathematicians on the basis of pure talent, G. H. Hardy \cite{H} gave Ramanujan the highest score of 100. On this scale, he gave Hilbert a score of 80, himself a score of only 25, but said his colleague Littlewood merited 30. In fact, Ramanujan's work has fascinated generations of mathematicians even a century after his death.
The curious identity recorded in one of his notebooks \begin{align*}
\sqrt{2 \left(1 - \frac{1}{3^2}\right)
\left(1 - \frac{1}{7^2}\right)
\left(1 - \frac{1}{11^2}\right)
\left(1 - \frac{1}{19^2}\right)}
\  =
\
 \left(1 + \frac{1}{7}\right) \left(1 + \frac{1}{11}\right) \left(1 + \frac{1}{19}\right)
\end{align*} was mentioned by Berndt \cite{Be} who asked: ``Is this an isolated result, or are there other identities of this type?'' Reb\'ak \cite{Re} provided formulas that generate infinitely many similar identities and believed that the curious identity is related to the reciprocal of the Landau-Ramanujan constant
\begin{align*}
\frac{1}{K}\ =\ \sqrt{2}\prod_{\substack{p\mbox{ prime}\\p\equiv 3\Mod 4}}\sqrt{1-\frac{1}{p^2}}\ =\ & \sqrt{2\left(1-\frac{1}{3^2}\right)\left(1-\frac{1}{7^2}\right)\left(1-\frac{1}{11^2}\right)\left(1-\frac{1}{19^2}\right)\cdots}.
\end{align*}
Chu \cite{Ch} examined necessary and sufficient conditions for the integers in the identity and proved that there are only finitely many such identities. A natural question to ask is whether we can twist the form of the original Ramanujan's identity to have infinitely many attractive identities. In this note, we answer this question positively. 

Given $t, h, x, y, z\in\mathbb{R}$, we define
\begin{align*}
    f(t, h, x, y, z) &\ =\ \sqrt{t \left(\frac{1}{h^2}-\frac{1}{x^2}\right)\left(\frac{1}{h^2}-\frac{1}{y^2}\right)\left(\frac{1}{h^2}-\frac{1}{z^2}\right)}, \mbox{ and }\\
    g(x,y,z)&\ =\ \left(\frac{1}{h}+\frac{1}{x}\right)\left(\frac{1}{h}+\frac{1}{y}\right)\left(\frac{1}{h}+\frac{1}{z}\right).
\end{align*}
We consider identities of the form
\begin{align}\label{mostgen}
f(t, h, x, y, z) \ =\ g(x,y,z).
\end{align}
Our first main result states necessary and sufficient conditions for real numbers $t, h, x, y$ and $z$ that satisfy Equation \eqref{mostgen}. Despite the complicated appearance of the equation, necessary and sufficient conditions for the variables can be established using quadratic equations. We then use the sufficient condition and choose values for the variables cleverly to have unexpected identities like the following
\begin{align*}
    \sqrt{8\left(\frac{1}{3^2}-\frac{1}{5^2}\right)\left(\frac{1}{3^2}-\frac{1}{15^2}\right)\left(\frac{1}{3^2}-\frac{1}{21^2}\right)}\ =\ \left(\frac{1}{3}+\frac{1}{5}\right)\left(\frac{1}{3}+\frac{1}{15}\right)\left(\frac{1}{3}+\frac{1}{21}\right)
\end{align*}
By a non-trivial solution to Equation \eqref{mostgen}, we require that
$$\begin{cases}
t, h, x, y, z\neq 0,\\
x, y, z\neq \pm h.
\end{cases}$$
\begin{theorem}\label{suffnec}
The tuple $(t, h, x, y, z)$ is a non-trivial solution to \eqref{mostgen} if and only if there exist $m, t, h, z\neq 0$ and $z\neq \pm h$ such that if we let 
\begin{align*}
    \gamma &\ =\ mt(z-h)-m(z+h),\\
    \beta &\ =\ hmt(z-h)+hm(z+h)-h^2,
\end{align*}
then 
    $$\begin{cases}\gamma^2-4\beta \ \ge\ 0,\\
   (h^2+h\gamma+\beta)\beta\left(\frac{1}{h}+\frac{1}{z}\right)\ >\ 0.\end{cases}$$
and $x$ and $y$ are the roots of $X^2-\gamma X+\beta = 0$.
\end{theorem}
%%%%%%%%%%%%%%%%%%%%%%%%%%%%%%%%%%%%%%%%%%%%%%%%%%%%%%%%%%%%%%%%%%%%%%%%%%%%%%%%%%%%%%%%%%%%%%%%%%%%%%%%%%%%%%%%%%%%%%%%%%%%%%%%%%%%%%%%%%%%%%%%%%%%%%%%%%%%%%%%%%%%%%%%%%%%%%%%%%%%%%%%%%%%%%%%%%%%%%%%%%%%%%%%%%%%%%%%%%%%%%%%%%%%%%%%%%%%%%%%%%%%%%%%%%%%%%%%%%%%%%%%%%%%%%%%%%%%%%%%%%%%%%%%%%%%%%%%%%%%
\section{Proof of Theorem \ref{suffnec}}
We first prove the forward implication. From $f(t,h,x,y,z) = g(x,y,z)$, we get 
$$t\left(\frac{1}{h}-\frac{1}{x}\right)\left(\frac{1}{h}-\frac{1}{y}\right)\left(\frac{1}{h}-\frac{1}{z}\right)\ =\ \left(\frac{1}{h}+\frac{1}{x}\right)\left(\frac{1}{h}+\frac{1}{y}\right)\left(\frac{1}{h}+\frac{1}{z}\right).$$
So, 
\begin{align*}
    t\frac{\left(\frac{1}{h}-\frac{1}{z}\right)}{\left(\frac{1}{h}+\frac{1}{z}\right)}\ =\ \frac{\left(\frac{1}{h}+\frac{1}{x}\right)\left(\frac{1}{h}+\frac{1}{y}\right)}{\left(\frac{1}{h}-\frac{1}{x}\right)\left(\frac{1}{h}-\frac{1}{y}\right)}\ =\ \frac{xy+h(x+y)+h^2}{xy-h(x+y)+h^2}.
\end{align*}
By simple arithmetic, we know that there exists $m\neq 0$ such that
\begin{align*}
    xy+h(x+y)+h^2 &\ =\ 2hmt(z-h)\\
    xy-h(x+y)+h^2&\ =\ 2hm(z+h)
\end{align*}
Hence, \begin{align*}xy &\ =\ hmt(z-h) + hm(z+h) - h^2 \ =\ \beta\\ x+y &\ =\ mt(z-h)-m(z+h)\ =\ \gamma.\end{align*}
Therefore, $x,y$ are the roots of $X^2-\gamma X+\beta = 0$, which implies that $\gamma^2-4\beta\ge 0$. Since a non-trivial solution makes both side of Equation \eqref{mostgen} strictly positive, $(h^2+h\gamma+\beta)\beta\left(\frac{1}{h}+\frac{1}{z}\right)>0$, as desired. 
%%%%%%%%%%%%%%%%%%%%%%%%%%%%%%%%%%%%%%%%%%%%%%%%%%%%%%%%%%%%%%%%%%%%%%%%%%%%%%%%%%%%%%%%%%%%%%%%%%%%

We now prove the backward implication. Since $x, y$ are the roots of $X^2-\gamma X + \beta$, we have
$x + y = \gamma$ and $xy  = \beta$. 
So,
\begin{align}
    \label{s1}xy + h(x+y) + h^2\ =\ 2hmt(z-h)\\
    \label{s2}xy - h(x+y) + h^2 \ =\ 2hm(z+h).
\end{align}
Divide \eqref{s1} by \eqref{s2} side by side to have
\begin{align*}
   \frac{t(z-h)}{z+h} \ =\ \frac{xy+h(x+y)+h^2}{xy - h(x+y) + h^2}\ =\ \frac{(x+h)(y+h)}{(x-h)(y-h)}.
\end{align*}
Because $x,y,h,z\neq 0$ (recall that $\beta\neq 0$), 
\begin{align*}t\frac{\frac{1}{h}-\frac{1}{z}}{\frac{1}{h}+\frac{1}{z}}\ =\ \frac{\left(\frac{1}{h}+\frac{1}{x}\right)\left(\frac{1}{h}+\frac{1}{y}\right)}{\left(\frac{1}{h}-\frac{1}{x}\right)\left(\frac{1}{h}-\frac{1}{y}\right)}.\end{align*}
Rearranging, we have
\begin{align}\label{hh}
    t\left(\frac{1}{h^2}-\frac{1}{x^2}\right)\left(\frac{1}{h^2}-\frac{1}{y^2}\right)\left(\frac{1}{h^2}-\frac{1}{z^2}\right)\ =\ \left(\left(\frac{1}{h}+\frac{1}{x}\right)\left(\frac{1}{h}+\frac{1}{y}\right)\left(\frac{1}{h}+\frac{1}{z}\right)\right)^2.
\end{align}
Since $(h^2+h\gamma+\beta)\beta\left(\frac{1}{h}+\frac{1}{z}\right)\ >\ 0$, we know that $\left(\frac{1}{h}+\frac{1}{x}\right)\left(\frac{1}{h}+\frac{1}{y}\right)\left(\frac{1}{h}+\frac{1}{z}\right)>0$ and so \eqref{hh} implies that $(t, h, x, y, z)$ is a non-trivial solution of \eqref{mostgen}.
%%%%%%%%%%%%%%%%%%%%%%%%%%%%%%%%%%%%%%%%%%%%%%%%%%%%%%%%%%%%%%%%%%%%%%%%%%%%%%%%%%%%%%%%%%%%%%%%%%%%%%%%%%%%%%%%%%%%%%%%%%%%%%%%%%%%%%%%%%%%%%%%%%%%%%%%%%%%%%%%%%%%%%%%%%%%%%%%%%%%%%%%%%%%%%%%%%%%%%%%%%%%%%%%%%%%%%%%%%%%%%%%%%%%%%%%%%%%%%%%%%%%%%%%%%%%%%%%%%%%%%%%%%%%%%%%%%%%%%%%%%%%%%%%%%%%%%%%%%%%
\section{Infinitely many Ramanujan-type identities}
We use the sufficient condition of Theorem \ref{suffnec} to give several forms of Ramanujan-type identities. Our goal is to choose suitable values of $m, t, h, z$ such that $\sqrt{\gamma^2-4\beta}$ is rational. We first give one example then give several unexpected formulas deduced by using the same method. 

Fix $h\in\mathbb{N}$. Let $t = 2+2h$, $z = h+2$ and $m = \frac{6h}{1+h}$. Then 
\begin{align*}
    \gamma &\ =\ mt(z-h)-m(z+h)\ =\ 12h,\\
    \beta &\ =\ hmt(z-h)+hm(z+h)-h^2\ =\ 35h^2.
\end{align*}
We have
\begin{align*}
    \gamma^2-4\beta \ =\ (12h)^2-4\cdot 35h^2 \ =\ 4h^2\ \ge\ 0,\\
    (h^2+h\gamma+\beta)\beta\left(\frac{1}{h}+\frac{1}{z}\right) \ >\ 0.
\end{align*}
Therefore, $m, t, h$ and $z$ satisfy all requirements of the sufficient condition. Because $x, y$ are the roots of $X^2-12hX+35h^2 = 0$, we have $(x,y) = (5h, 7h)$ and we get the following identity
\begin{align}\label{u1}
    &\sqrt{(2+2h)\left(\frac{1}{h^2}-\frac{1}{(h+2)^2}\right)\left(\frac{1}{h^2}-\frac{1}{(5h)^2}\right)\left(\frac{1}{h^2}-\frac{1}{(7h)^2}\right)}\nonumber\\
    &\ =\ \left(\frac{1}{h}+\frac{1}{h+2}\right)\left(\frac{1}{h}+\frac{1}{5h}\right)\left(\frac{1}{h}+\frac{1}{7h}\right).
\end{align}
Hence, we have infinitely many Ramanujan-type identities with different values of $h\in\mathbb{N}$. The following table provides some other values of $m, t, h, x, y, z$ that give \eqref{mostgen}. These are found by trial and error with the help of computers. 
%%%%%%%%%%%%%%%%%%%%%%%%%%%%%%%%%%%%%%%%%%%%%%%%%%%%%%%%%%%%%%%%%%%%%%%%%%%%%%%%%%%%%%%%%%%%%%%%%%%%
\begin{center}
\begin{tabular}{ |c|c|c|c|c| }
\hline
$t$ & $m$ & $x$ & $y$ & $z$ \\
\hline
$2+2h$ & $6h/(h+1)$ & $5h$ & $7h$ & $h+2$ \\
 $2+2h$ & $15h/(2h+2)$ & $4h$ & $11h$ & $h+2$ \\
 $2+2h$ & $28h/(3h+3)$ & $11h/3$ & $15h$ & $h+2$ \\
 $2+2h$ & $35h/(6h+6)$ & $17h/3$ & $6h$ & $h+2$\\
 $2+2h$ & $45h/(4h+4)$ & $7h/2$ & $19h$ & $h+2$ \\
 $2+2h$ & $20h/(3h+3)$ & $13h/3$ & $9h$ & $h+2$ \\
 \hline
 $3+2h$ & $4h/(2h+3)$ & $3h$ & $5h$ & $h+3$ \\
 $3+2h$ & $4h/(2h+3)$ & $7h/2$ & $4h$ & $h+3$ \\
 \hline
\end{tabular}
\end{center}
\begin{remark}
The next-to-last row gives 
\begin{align}
    &\label{u2} \sqrt{(3+2h)\left(\frac{1}{h^2}-\frac{1}{(h+3)^2}\right)\left(\frac{1}{h^2}-\frac{1}{(3h)^2}\right)\left(\frac{1}{h^2}-\frac{1}{(5h)^2}\right)}\nonumber\\
    &\ =\ \left(\frac{1}{h}+\frac{1}{h+3}\right)\left(\frac{1}{h}+\frac{1}{3h}\right)\left(\frac{1}{h}+\frac{1}{5h}\right).
\end{align}
Furthermore, a quick computer search gives 
\begin{align}\label{u3}
    \sqrt{3\left(\frac{1}{2^2}-\frac{1}{5^2}\right)\left(\frac{1}{2^2}-\frac{1}{22^2}\right)\left(\frac{1}{2^2}-\frac{1}{58^2}\right)} \ =\ \left(\frac{1}{2}+\frac{1}{5}\right)\left(\frac{1}{2}+\frac{1}{22}\right)\left(\frac{1}{2}+\frac{1}{58}\right).
\end{align}
Identities \eqref{u1}, \eqref{u2} and \eqref{u3} show that for all $t\ge 3$, there exists a non-trivial, positive, integral solution $(t, h, x, y, z)$ to \eqref{mostgen}. It is not hard to show that when $t = 1$ or $2$, \eqref{mostgen} does not have any non-trivial, positive, integral solution. 
\end{remark}

%%%%%%%%%%%%%%%%%%%%%%%%%%%%%%%%%%%%%%%%%%%%%%%%%%%%%%%%%%%%%%%%%%%%%%%%%%%%%%%%%%%%%%%%%%%%%%%%%%%%%%%%%%%%%%%%%%%%%%%%%%%%%%%%%%%%%%%%%%%%%%%%%%%%%%%%%%%%%%%%%%%%%%%%%%%%%%%%%%%%%%%%%%%%%%%%%%%%%%%%%%%%%%%%%%%%%%%%%%%%%%%%%%%%%%%%%%%%%%%%%%%%%%%%%%%%%%%%%%%%%%%%%%%%%%%%%%%%%%%%%%%%%%%%%%%%%%%%%%%%

\bigskip
\hrule
\bigskip

\noindent 2010 {\it Mathematics Subject Classification}: Primary 11A67; Secondary 11E25.

\noindent \emph{Keywords:} Ramanujan identity, reciprocal, square root

\bigskip
\hrule
\bigskip

\end{document}